\documentclass[12pt]{amsart}

\usepackage{amsmath,amsfonts,amssymb,amscd}
\begin{document}
\renewcommand{\thefootnote}{\fnsymbol{footnote}}
\pagestyle{plain}
\title{Extremal metrics and stabilities\\
 on polarized manifolds}
\author{Toshiki Mabuchi${}^*$}
\address{Department of Mathematics, Graduate School of Science, Osaka University, 
Toyonaka, Osaka, 560-0043 Japan
}
\maketitle
\footnotetext{Special thanks are due to Professors K. Hirachi and G. Komatsu for valuable suggestions. }
\section{Introduction} 

Let $M$ be a compact complex connected manifold.
As an introduction to our subjects,
we recall the following well-known conjecture of Calabi \cite{Cal1}:

\medskip\noindent
{\bf Conjecture}.
{\rm (i)} 
{\em If $\,c_1 (M)_{\Bbb R} < 0$, then $M$ admits a unique K\"ahler-Einstein metric $\omega$ 
such that $\operatorname{Ric}(\omega ) = - \omega$.} \newline
{\rm (ii)} 
{\em If $\,c_1 (M)_{\Bbb R} = 0$, then each K\"ahler class on $M$ admits a unique
K\"ahler-Einstein metric
$\omega$  such that $\operatorname{Ric}(\omega ) = 0$.}
\newline
{\rm (iii)} {\em For $\,c_1 (M)_{\Bbb R} > 0$, find a suitable condition for $M$ to admit 
a K\"ahler-Einstein metric $\omega$ 
such that $\operatorname{Ric}(\omega ) =  \omega$.}

\medskip
Based on the pioneering works of Calabi \cite{Cal2}, \cite{Cal3} and Aubin \cite{Aub}, 
a complete affirmative answer to (i) and (ii) was given by Yau \cite{Yau} 
by solving systematically
certain complex Monge-Amp\`ere equations. However for (iii), sufficient conditions are 
known only partially
by Siu \cite{Siu}, Nadel \cite{Nad}, Tian \cite{Tia1}, \cite{Tia}, Tian and Yau \cite{TY}, 
Wang and Zhu \cite{WZ}, where some necessary conditions were formulated as obstructions by Futaki \cite{Fut1} and Matsushima \cite{Mat} (see also Lichnerrowicz \cite{Lic}). This (iii) mainly motivates our studies of 
stabilities and extremal metrics, while
an analysis of destabilizing phenomena
 caused by nonexistence of K\"ahler-Einstein metrics allows us to obtain
very nice by-products such as Nadel's vanishing theorem  \cite{Nad}
via the use of multiplier ideal sheaves.
 
 \section{Stability for manifolds in algebraic geometry}
 
 In Mumford's GIT \cite{MFK}, 
 moduli spaces of algebraic varieties are constructed
via the theory of invariants, where varieties are described by numerical 
 data modulo actions of reductive algebraic groups.
Then roughly speaking, stable points are those away from very bad points in moduli spaces. 
For a precise definition, consider a representation of 
 a reductive algebraic group $G$
 on a complex vector space $W$.
Let $w \in W$. 
We denote by $G_w$ the isotropy subgroup
 of $G$ at $w$.

\medskip\noindent 
{\bf Definition 2.1.}
A point $w \in W$ 
 is said to be {\it stable\/} (resp. {\it properly stable\/} )
 if the orbit $G \cdot w$ is closed in $W$ 
 (resp. $G \cdot w$ is closed in $W$ with 
 $|G_w| < \infty$).
\medskip
 
 For moduli spaces of polarized varieties, 
the Chow-Mumford stability and the Hilbert-Mumford stability are known.
Hereafter, by a polarized manifold $(M, L)$, we mean a very ample holomorphic line bundle $L$ over 
a nonsingular projective algebraic variety $M$ defined over $\Bbb C$,
where arguments in this section has nothing to do with the nonsingularity of $M$. 
Now for a polarized manifold $(M,L)$, put $n:= \dim M$
and let $m$ be a positive integer.
Then associated to the complete linear system $|L^m|$,
we have the Kodaira embedding 
$$
\iota_m  : M\; \hookrightarrow \;\Bbb P^*(V_m),
$$
where $\Bbb P^*(V_m)$ denotes the set of all hyperplanes in 
$V_m := H^0(M, \mathcal{O}(L^m))$ through the origin.
Let $d_m$ be the
degree of  $\iota_m (M )$ in the projective space $\Bbb P^* (V_m)$. 
Put $G_m := \operatorname{SL}_{\Bbb C}(V_m)$ 
and $W_m :=\{S^{d_m} (V_m )\}^{\otimes n+1}$, 
where $S^{d_m} (V_m )$ 
is the 
$d_m$-th symmetric tensor product of the space $V_m$.
Take an element
$M_m \neq 0$ in $W_m^*$ such that 
the associated element $[M_m ]$ in 
$\Bbb P^* (W_m)$ is the Chow point of the 
irreducible reduced algebraic cycle
$\iota_m (M)$ on $\Bbb P^* (V_m )$. 
For the natural action of $G_m$ on $W_m^*$,
we now apply Definition 2.1 above to $G = G_m$ and $W = W_m$:

\medskip\noindent
{\bf Definition 2.2.}
(a)\, $(M,L^m)$ is said to be {\it Chow-Mumford stable\/} (resp. 
{\it Chow-Mumford properly stable\/})
if $M_m$ in $W_m^*$ is stable (resp. properly stable).

\smallskip\noindent
(b)\, $(M, L)$ is said to be {\it asymptotically Chow-Mumford stable\/}
(resp. {\it asymptotically Chow-Mumford properly stable\/})
if, for $m \gg 1$, 
$(M, L^m)$ is Chow-Mumford stable (resp. Chow-Mumford properly stable).

\medskip
Let $m$ and $k$ be positive integers.
Then the kernel $I_{m,k}$ of the natural homomorphism 
of $S^k (V_m )$ to $V_{m k} := H^0(M, \mathcal{O}(L^{m k}))$
is the homogeneous ideal of degree $k$ defining $M$ in 
$\Bbb P^*(V_m)$. Put $N_{mk} := \dim V_{m k}$ and $\gamma_{m,k} := \dim I_{m,k}$.
Then $\wedge^{\gamma_{m,k}} I_{m,k}$ is a complex line 
in the vector space $W_{m,k}:=\wedge^{\gamma_{m,k}} (S^k (V_m ))$.
Take an element $w_{m,k}\neq 0$ in $\wedge^{\gamma_{m,k}} I_{m,k}$. 
For the natural action of $G_m := \operatorname{SL}_{\Bbb C} (V_m)$ on $W_{m, k}$, 
we apply Definition 2.1 to $G = G_m$ and $W = W_{m, k}$: 

\medskip\noindent
{\bf Definition 2.3.} 
(a)\, $(M,L^m)$ is said to be {\it Hilbert-Mumford stable\/} (resp.
{\it Hilbert-Mumford properly stable\/}) if
$w_{\ell,k}\in W_{\ell,k}$ is stable (resp. properly stable).

\smallskip\noindent
 (b)\, $(M,L)$ is said to be {\it asymptotically 
Hilbert-Mumford stable} (resp. {\it asymptotically 
Hilbert-Mumford properly stable}) if, 
for all $m \gg 1$,  
 $(M,L^m)$ is Hilbert-Muford stable (resp. Hilbert-Mumford properly stable).

\medskip
A result of Fogarty \cite{Fog} shows that, if $(M, L^m)$ 
is Chow-Mumford stable, 
then $(M, L^m)$ is also 
Hilbert-Mumford stable. 
Though the converse has been unknown,
the relationship between these two stabilities
are now becoming clear
 (cf. \cite{Mab6}).

\medskip
Stability for manifolds is an important subject
in moduli theories of algebraic geometry.
Recall, for instance, the following famous result 
of Mumford \cite{Mum}: 

\medskip\noindent
{\bf Fact 1}: {\em If $\, L$ is an ample line bundle of degree 
$d \geq 2 g +1$ 
over a compact Riemann surface $C$ of genus $g \geq 1$, then $\, (C, L)$ is Chow-Mumford properly stable.}

\medskip
For the pluri-canonical bundles $K^{\otimes m}_M$ on $M$, $m \gg 1$, 
using the asymptotic Hilbert-Mumford stability, Gieseker \cite{Gie} generalized this result to the case where $M$ is a surface of general type.
For higher dimensions, a stability result by Viehweg \cite{Vie} is known in the case where the canonical bundle $K_M$ is semipositive.  However, both for the
results of Gieseker and Viehweg, the proof of stability is fairly complicated, while the underlying manifold (or orbifold) admit a K\"ahler-Einstein metric.

\section{The Hitchin-Kobayashi correspondence and its manifold analogue}

\smallskip
For a holomorphic vector bundle $E$ over an $n$-dimensional
 compact K\"ahler manifold $(M, \omega )$,
we say that $E$ is {\it Takemoto-Mumford stable\/} if 
$$
\frac{\int_M c_1(\mathcal{S})\, \omega^{n-1}}{\operatorname{rk}(\mathcal{S})}
\; <\; \frac{\int_M c_1 (E) \,\omega^{n-1}}{\operatorname{rk}(E)}
$$
for every coherent subsheaf $\mathcal{S}$ of $\mathcal{O}(E)$ satisfying 
$0 < \operatorname{rk}(\mathcal{S}) < \operatorname{rk}(E)$.
Recall the following Hitchin-Kobayashi correspondence for 
vector bundles:

\medskip\noindent
{\bf Fact 2}:
{\em An indecomposable holomorphic vector bundle $E$ over $M$
is Takemoto-Mumford stable if and only if $E$ admits a Hermitian-Einstein
metric.}

\medskip
This fact was established in 1980's by Donaldson \cite{Don1}, 
Kobayashi \cite{Kob}, L\"ubke \cite{Lub}, Uhlenbeck and Yau \cite{UY}. 
As a manifold analogue of this conjecture, we can naturally ask whether the following
conjecture  (known as Yau's conjecture)  is true:

\medskip\noindent
{\bf Conjecture}. 
{\em The polarization class of $\, (M, L)$
admits a K\"ahler metric of constant scalar curvature $\,($or more generally an extremal K\"ahler 
metric$)$ if and only if $\,(M, L)$ is asymptotically stable in a certain sense of GIT.}

\medskip
For ``only if'' part of this conjecture, the first breakthrough 
was made by Tian \cite{Tia3}. By introducing the concept of K-stability,
he gave an answer to the ``only if'' part for K\"ahler-Einstein manifolds,
and showed that
some Fano manifolds without nontrivial holomorphic vector fields
admit no K\"ahler-Einstein metrics.
A remarkable progress was made by Donaldson \cite{Don2}
who showed the Chow-Mumford stability for a polarized 
K\"ahker manifold $(M, \omega )$ of constant scalar curvature essentially when the connected linear algebraic part $H$
of  the group $\operatorname{Aut}(M)$ of holomorphic automorphisms of $M$ is
semisimple.  In the present paper, we shall show how Donaldson's work is generalized to extremal K\"ahler cases  without any assumption
of $H$ (see also \cite{Mab4}, \cite{Mab5}). The relationship between this 
generalization and  a recent result by Chen and Tian \cite{CT2} will be treated 
elsewhere.

\section{The asymptotic Bergman kernel}

For a polarized manifold $(M, L)$, take a Hermitian metric $h$ for $L$ such that 
$\omega := c_1 (L; h)$ is a K\"ahler form.  Define a Hermitian pairing on 
$V_m := H^0(M, \mathcal{O}(L^m))$ by
$$
<\sigma_1, \sigma_2 >^{}_{L^2}\; := \; \int_M (\sigma_1, \sigma_2)_h \,\omega^n,
\qquad \sigma_1, \sigma_2  \in V_m,
$$
where $(\; , \; )_h$ denotes the pointwise Hermitian pairing by $h^m$
for sections for $L^m$. For an orthonormal basis $\{ \sigma_1, \sigma_2, \dots, \sigma_{N_m}\}$ 
of $V_m$, we put
\begin{equation}
B_{m,\omega} \;:=\; \frac{n!}{\; m^n}\, (|\sigma_1|_h^2 + |\sigma_2 |_h^2 + \dots + |\sigma_{N_m} |_h^2),
\end{equation}
where $|\sigma |_h^2 := (\sigma,\sigma )_h$ for $\sigma \in V_m$. This $B_{m,\omega}$ is called the $m$-th 
{\it Bergman kernel\/} for $(M, \omega )$ (cf. Tian \cite{Tia2}, Zelditch \cite{Zel}, Catlin \cite{Cat}),
where we consider the asymptotic behavior of $B_{m,\omega}$
as $m \to \infty$.
Note that $B_{m,\omega}$ depends only on $(m, \omega )$,
and is independent of the choice of both $h$ and the orthonormal basis for $V_m$. Next for
$$
D \, :=\, \{\,\ell \in L^*\,;\, |\ell |_h < 1\,\},
$$
the boundary $X := \partial D = \{\ell \in L^*; |\ell |_h = 1\}$ over $M$ is an $S^1$-bundle.
Let $\operatorname{pr} : X \to M$ be the natural projection.
We now consider the Szeg\"o kernel 
$$
S_{\omega} := S_{\omega}(x,y)
$$ 
for the 
projection of $L^2(X)$ onto the Hardy space $L^2(X) \cap \Gamma (D, \mathcal{O})$ of  
boundary values of holomorphic functions on $D$. 
Then for each positive integer $m$, 
the corresponding $m$-th Bergman 
Kernel $B_{m,\omega}$ for the K\"ahler manifold $(M,\omega )$ is characterized as the Fourier coefficient
$$
\operatorname{pr}^*B_{m,\omega } \; :=\; \; \frac{n!}{\; m^n}\,
 \int_{S^1} e^{-i m \theta}\;
S_{\omega}(e^{i \theta }x, x)\; d\theta.
$$
Now the Bergman kernel is defined not only for positive integers $m$, 
but also for complex numbers $\xi$ as follows. To see the situation, we first 
consider the case where $M$ is a single point. Then in place of $S_{\omega}(e^{i \theta }x, x)$,  consider a smooth function $S = S(\theta )$ on $S^1:= \Bbb R /2 \pi \Bbb Z$ for simplicity. The associated Fourier coefficient $B_m$ is 
$$
B_m \;=\; \int_{S^1} e^{-i m \theta}\;
S (\theta )\; d\theta
$$
for each integer $m\neq 0$.
Then for open intervals $I_1 := (-3\pi/4 ,  3\pi/4 )$ and $I_2 = (\pi/4, 7\pi/4 )$ in $\Bbb R$, 
we choose the open cover
$$
S^1 = U_1 \cup U_2
$$ 
where $U_1 := I_1$ mod $2\pi$ and $U_2 := I_2$ mod $2\pi$. 
By choosing a partition of unity subordinate to this
open cover, we write 
$$
\rho_1 (\theta )  + \rho_2 (\theta ) \; =\; 1,
\qquad \theta \in S^1,
$$
where $\rho_{\alpha}\in C^{\infty}(S^1)_{\Bbb R}$, $\alpha = 1,2$, 
are functions $\geq 0$ satisfying
$\operatorname{Supp}(\rho_\alpha ) \subset U_\alpha$.
For the coordinate $\tilde{\theta}$ for $\Bbb R$, writing
$\tilde{\theta}$ mod $2\pi$ as $\theta$, define
$\tilde{\rho}_{\alpha}\in C^{\infty}(\Bbb R )_{\Bbb R}$, $\alpha = 1,2$, by
$$
\tilde{\rho}_{\alpha} (\tilde{\theta} ) \; =\; 
\left \{ \begin{array}{ll}  \;\; \rho_{\alpha} (\theta ),   \qquad & \tilde{\theta} \in I_{\alpha} \\
\;\; \;\; 0, \qquad &\tilde{\theta} \notin I_{\alpha}.
\end{array}
\right.
$$
Then the Fourier transform $\mathcal{F}(S) = \mathcal{F}(S)(\xi )$
of $S$ is an entire function in $\xi\in \Bbb C$ defined as the integral
$$
\mathcal{F}(S) (\xi ) \; =\; \int_{\Bbb R} e^{-i \xi \theta}\;
\{\, \rho_1 (\tilde{\theta}) + \rho_2 (\tilde{\theta}) \,\} \;
S (\tilde{\theta} )\; d \tilde{\theta},
\qquad \xi \in \Bbb C,
$$
satisfying $\mathcal{F}(S) (m ) = B_m$ for all integers $m$. 
Though $\mathcal{F}(S)$ may depends on the choice of the 
partition of unity, its restriction to $\Bbb Z$ is unique.
This situation is easily understood for instance by the fact the functions $\mathcal{F}(S) (\xi )$ 
and $\mathcal{F}(S) (\xi ) + \sin (\pi \xi )$ in $\xi$ coincide 
on $\Bbb Z$.

\medskip
Now the above process of generalization from the Fourier series to the Fourier transform 
 is valid also for the case where the base space $M$ is nontrivial. Actually, we define an
entire function 
$B_{\xi, \omega}$ in $\xi\in \Bbb C$ by
$$
\operatorname{pr}^* B_{\xi,\omega } \; :=\; \; \frac{n!}{\; \xi^n}\,
 \int_{\Bbb R} \; e^{-i \xi \tilde{\theta}}\;
 \{\, \rho_1 (\tilde{\theta}) + \rho_2 (\tilde{\theta}) \,\} \;
S_{\omega}(e^{i \tilde{\theta} }x, x)\; d\tilde{\theta}.
$$
By setting $q := \xi^{-1}$,  we study the asymptotic behavior of 
$B_{\xi, \omega}$, as $q\to 0$ along the positive real line
$\{\, q >0 \,\}$.
Let $\sigma_{\omega}$ denote the scalar curvature of the K\"aler manifold $(M, \omega)$.
Then as in discrete cases by Tian \cite{Tia2}, Zelditch \cite{Zel}, Catlin \cite{Cat}, the asymptotic expansion of 
$B_{\xi, \omega}$ in $q$ yields
\begin{equation}
B_{\xi,\omega } \; =\;  1 + a_1 (\omega) q + a_2(\omega) q^2 + \dots,
\qquad 0\leq  q \ll 1,
\end{equation}
where $a_1(\omega) = \sigma_{\omega}/2$ by a result of Lu \cite{Lu}. 
For more details of the expansion in discrete cases, see also Hirachi \cite{Hir}.

\section{Balanced metrics}

Choose a Hermitian metric $h$ for $L$ such that
$\omega := c_1(L; h)$ is a K\"ahler form.
Then $\omega$ is called an $m$-th {\it balanced metric\/} (cf. \cite{Zha}, \cite{Luo})\, 
for $(M,L)$, if $B_{m,\omega}$ is a constant function $( = C_m )$ on $M$.
First, put $q := 1/m$. By integrating (1) and (2) on the 
K\"ahler manifold $(M, \omega )$, we see that $C_m$ is written as
$$
C_m := \frac{n!}{m^n c_1(L)^n[M]}\, N_m =   1 + \frac{n}{2}\frac{c_1(M) c_1(L)^{n-1}[M]}{c_1(L)^n[M]}\,q\, + 
\, O(q^2), \quad 0 \leq q \ll 1,
$$
where the left-hand side is the Hilbert polynomial $P(m)$ for $(M,L)$ 
divided by $m^n c_1(L)^n[M]/n!$. Hence, it is easy to define $C_{\xi}$ by setting
\begin{equation}
C_{\xi} \; =\;\frac{ n! \, P(\xi ) }{\xi^n c_1(L)^n[M]}, 
\quad \qquad \xi \in \Bbb C^*.
\end{equation}
Put $\Delta_{\omega} := \bar{\partial}^*\bar{\partial}+ \bar{\partial}\bar{\partial}^*$ on $(M,\omega)$. Define the modified Bergman kernel $\beta_{q,\omega}$ by
\begin{equation}
\beta_{q,\omega} \; :=\; 2\xi \,(1+\frac{2q}{3}\,\Delta_{\omega})\, (B_{\xi,\omega} - C_{\xi} )
\; =\; \sigma_\omega - \bar{\sigma}_\omega + O(q),
\end{equation}
where the average $\bar{\sigma}_{\omega}$ of the scalar curvature $\sigma_{\omega}$ is $n c_1(M) c_1(L)^{n-1}[M]/c_1(L)^n[M]$ independent of the choice of $\omega$ in $c_1(L)_{\Bbb R}$, and we now put $ q := 1/\xi$.
Then for $\xi =m$, $\omega$ is an $m$-th balanced metric for $(M, L)$ if and only if $\beta_{q,\omega}$  vanishes everywhere on $M$. 
Recall the following result by Zhang \cite{Zha} (cf. \cite{Luo}; see also 
\cite{Mab3}):

\medskip\noindent
{\bf Fact 3}: 
{\em $(M,L^m)$ is Chow-Mumford stable if and only if $(M, L)$ admits an $m$-th 
balanced metric.}

\medskip\noindent
Consider the maximal connected linear algebraic subgroup $H$ of $\operatorname{Aut}(M)$,
so that the identity component of $\operatorname{Aut}(M)/H$ is an Abelian variety.
Let us now choose an algebraic torus $T \cong (\Bbb C^*)^r$  
in the connected component $Z^{\Bbb C}$ of the center 
of the reductive part $R(H)$ for $H$ 
in the Chevalley decomposition 
$$
H \; =\; R(H)  \ltimes  H_u,
$$
where $H_u$ is the unipotent 
radical of $H$. Replacing $L$ by its suitable positive multiple, we may asssume that the $H$-action 
on $M$ is lifted to a bundle action on $L$ covering the $H$-action on $M$.
For each character $\chi \in \operatorname{Hom}(T, \Bbb C^*)$, we set
$$
W \chi  \; :=\; \{\, \sigma \in V_m\,;\,  \sigma\cdot g = \chi (g) \sigma
\; \text{ for all $g \in T$}\,\},
$$
where $V_m \times T \owns (\sigma , g) \mapsto  \sigma\cdot g$ is the right $T$-action
on $V_m = H^0(M, \mathcal{O}(L^m))$
induced by the left
$H$-action on $L$. 
Now we have characters $\chi_k
\in \operatorname{Hom}(T, \Bbb C^*)$, $k =1,2,\dots, r_m$, such that the vector space $V_m$ is expressible
as a direct sum
$$
V_m \;=\; \bigoplus_{k=1}^{\; r_m} \; W\chi_k.
$$
For the maximal compact torus $T_c$ in $T$, 
we may assume that both $h$ and $\omega$ are $T_c$-invariant.
Put $\Bbb J_m := \{1,2,\dots, N_m\}$, where $N_m := \dim V_m$.
Choose an orthonormal basis $\{\,\sigma_1,\sigma_2,
\dots,\sigma_{N_m}\,\}$ for $V_m$ such that
all $\sigma_j$, $j\in \Bbb J_m$, belong to 
the union $\cup_{k=1}^{r_m} W\chi_k $.
 Hence, there exists a map
$\kappa : \Bbb J_m   \to  \{1,2,\dots, r_m\}$
satisfying
$$
\sigma_j \in W\chi_{\kappa (j)},
\qquad j \in \Bbb J_m.
$$  
Put $\frak t_{\Bbb R} := \,i\,\frak t_c$ for the Lie algebra $\frak t_c$ of $T_c$,
where $i := \sqrt{-1}$.
For each $\mathcal{Y}\in \frak t_{\Bbb R}$, 
by setting $g: = \exp (\mathcal{Y}/2)$, we put 
$h_g := h\cdot g$ for the natural $T$-action on the space of Hermitian metrics on $L$.
Define the $m$-th weighted Bergman kernel $B_{m,\omega, \mathcal{Y}}$, 
twisted by $\mathcal{Y}$, for 
$(M,\omega )$ by setting
$$
B_{m,\omega, \mathcal{Y}}\; :=\; \frac{n!}{\; m^n} \; 
\sum_{j=1}^{N_m} |\sigma_j |_{h_g}^2\; 
=\; g^*\left\{ \frac{n!}{\; m^n} \; 
\sum_{j=1}^{N_m}\; 
\frac{\; |\sigma_j |_h^2}{|\chi^{}_{\kappa (j)} (g) |^2}\right\}.
$$
Then $\omega$ is called an {\it $m$-th $T$-balanced metric\/} on $(M, L)$ if 
$B_{m,\omega, \mathcal{Y}}$ is a constant function (= $C_{m,\mathcal{Y}}$) on $M$ for some $\mathcal{Y}\in \frak t$. Consider the natural action of the group
$$
G_m\;  := \; \bigoplus_{k=1}^{r_m}\operatorname{SL}_{\Bbb C}(W\chi_k )
$$
acting on $V_m = \oplus_{k=1}^{r_m} W\chi_k$ diagonally (factor by factor). We say that $(M, L^m)$ is   
{\it Chow-Mumford $T$-stable\/} if the orbit 
$G_m \cdot M_m$ is closed in $W_m^*$.  
Note that (cf. \cite{Mab3})

\medskip\noindent
{\bf Fact 4}: If $M$ admits an $m$-th $T$-balanced metric on $(M,L)$, then
$(M,L^m)$ is Chow-Mumford $T$-stable.

\medskip
We shall now extend $\{\,B_{m, \omega, \mathcal{Y}}\,;\, m = 1,2,\dots \}$ to 
$\{\, B_{\xi, \omega, \mathcal{Y}}\,;\, \xi \in \Bbb C\,\}$
in such a way that 
$B_{m, \omega, \mathcal{Y}}$ coincides with ${B_{\xi, \omega, \mathcal{Y}}}_{| \xi = m}$ for all positive integers $m$. By the definition of $B_{m, \omega, \mathcal{Y}}$ (see also \cite{Mab3}, p.578) the equality
$B_{m,\omega, \mathcal{Y}} = B_{m,\omega} (h_g/h)^m$ always holds. Hence we put
\begin{equation}
B_{\xi, \omega,\mathcal{Y}} \; := \; B_{\xi, \omega} \cdot (h_g/h)^{\xi}\;
=\; B_{\xi, \omega}  \exp\{ \xi \log (h_g/h)\},
\qquad \xi \in \Bbb C.
\end{equation}
Once $B_{\xi, \omega,\mathcal{Y}}$ is defined, we can also define 
$C_{\xi,\mathcal{Y}}$, $\xi \in \Bbb C$, in such a way that 
$C_{m,\mathcal{Y}}$ coincides with ${C_{\xi, \mathcal{Y}}}_{|\xi = m} $.
Actually, we put
\begin{equation}
C_{\xi, \mathcal{Y}} := 
\int_M \; \frac{B_{\xi, \omega,\mathcal{Y}}}{c_1(L)^n[M]} \; g^*\omega^n.
\end{equation}

\section{A simple heuristic proof of Donaldson's theorem}

In this section, we shall show that  a heuristic application of the implicit function theorem simplifies the proof of Donaldson's theorem \cite{Don2}.
Fix a K\"ahler metric $\omega_0$ in $c_1(L)_{\Bbb R}$ of constant scalar curvature.
Assume that the group $H$ in the previous section is trivial.
For each K\"ahler metric $\omega$ in $c_1(L)_{\Bbb R}$, we can associate 
a unique real-valued smooth function $\varphi$ on $M$ such that
$$
\omega\; =\; \omega_0 + \sqrt{-1}\partial\bar{\partial}\varphi
$$
with normalization condition $\int_M \varphi\, \omega_0^n = 0$.
For an arbitrary nonnegative integer $k$ and a real number 
$\alpha$ satisfying $0<\alpha <1$, we more generally 
consider the case where $\varphi \in C^{k+4, \alpha}(M)_{\Bbb R}$, so that 
$\omega$ is a $C^{k+2, \alpha}$ K\"ahler metric on $M$.
The Fr\'echet derivative $D_\omega \sigma_{\omega}$ 
at $\omega = \omega_0$ of the scalar curvature function 
$\omega \mapsto \sigma_{\omega}$ is given by
$$
\left \{\,({D^{}_{\omega}\sigma^{}_{\omega}})\,
(\sqrt{-1}\,\partial\bar{\partial} \varphi )\,\right\}_{| \omega = \omega_0}
\; =\; \lim_{\varepsilon \to 0}
\frac{\sigma_{\omega_0 + \sqrt{-1}\varepsilon \partial
\bar{\partial}\varphi}  \, -\, \sigma_{\omega_0}}{\varepsilon}
\;=\; L_{\omega_0} \varphi,
$$
where $L_{\omega_0} : C^{k+4, \alpha}(M)_{\Bbb R} \to C^{k, \alpha}(M)_{\Bbb R}$ is the Lichnerowicz operator for the 
K\"ahler metric $\omega_0$ (cf. \cite{Lic}, \cite{Don2}).
Then by (4), the Fr\'echet derivative $D_\omega \beta_{q, \omega}$ of $\beta_{q, \omega}$ with respect to $\omega$
at $(q, \omega) = (0, \omega_0 )$ is written as
$$
\left \{\, (D_\omega \beta_{q, \omega})\,(\sqrt{-1}\,\partial\bar{\partial}\varphi )\,\right \}_{|(q,\omega ) = (0, \omega_0 )} \; =\; 
\left \{\,({D^{}_{\omega}\sigma^{}_{\omega}})\,
(\sqrt{-1}\,\partial\bar{\partial} \varphi )\,\right\}_{| \omega = \omega_0}
\; =\; L_{\omega_0} \varphi,
$$
where $L_{\omega_0}$ is an invertible operator by the triviality of $H$. 
By setting $q := 1/\xi$, we move $q$ in the half line 
$\Bbb R_{\geq 0} := \{0\} \cup \{\,1/\xi\,;\,\xi >0 \,\}$.
Replacing $q$ by $q^2$ if necessary,  
we apply the implicit function theorem to the map 
$(q, \omega )  \mapsto  \beta_{q,\omega}$.
(The required regularity for this map is rather delicate: By using \cite{BS}, Theorem 1.5 and \S 2.c,
we can write both $B_{\xi,\omega}$ and its $\omega$-derivative 
as integrals similar to (18) in \cite{Zel}. Then
the estimate of the remainder term in the asymptotic expansion 
for $B_{m,\omega}$ in \cite{Zel}, Theorem 1, 
is valid also for the asymptotic expansion of $B_{\xi,\omega}$ and its $\omega$-derivative.
However, we need more delicate estimates to see the continuity for
$\beta_{q,\omega}$ and its $\omega$-derivative. Related to Nash-Moser's process,
this will be treated elsewhere.)
Then we have openness of the solutions
for the one-parameter family of equations
\begin{equation}
\beta_{q,\omega} \; =\; 0,
\qquad q \geq 0,
\end{equation}
i.e., there exists a one-parameter family of $C^{k+2,\alpha}$ solutions 
$\omega = \omega (q)$, $0\leq q <\varepsilon$, for (7) with
$\omega (0) = \omega_0$,  where 
$\varepsilon$ is sufficiently small.
Hence by Fact 3 in Section 5, $(M,L^m)$ is Chow-Mumford 
stable for all integers $m > 1/\varepsilon$.
\qed

\section{The case where $M$ admits symmetries}

In this section, we consider 
a polarized manifold $(M,L)$ with an extremal K\"ahler metric 
$\omega_0$ in $c_1 (L)_{\Bbb R}$. Then following Section 6,
 a result on stability
in \cite{Mab5} will be discussed.
Let $\mathcal{V}$ be the 
associated extremal K\"ahler vector field on $M$.
Assume that $H$ 
is possibly of positive dimension. 
 Then the identity component $K$ 
of the isometry group 
of $(M, \omega_0)$ is a maximal compact subgroup 
in $H$. 
Let $\frak z$ be the Lie algebra of the 
of the identity component $Z$ of the center of $K$.

\medskip\noindent
{\em Step 1}.  As the first step, we assume that
the algebraic torus $T$ in Section 5 
is the complexification $Z^{\Bbb  C}$ of $Z$ in $H$,
so that the Lie algebra $\frak t$ of $T$ coincides with
the Lie algebra $\frak z^{\Bbb C}$ of $Z^{\Bbb  C}$.
Let $q\in \Bbb R_{\geq 0}$, where
we set $q = 1/\xi$ for the part $0<\xi \in \Bbb R$.
Take an element $\mathcal{W}$ in the Lie algebra $\frak{z}$.
Let $\omega$ be a K\"ahler metric in the
class $c_1(L)_{\Bbb R}$. 
Then by setting $\mathcal{Y}:= \,i\, q^2(\mathcal{V} + \mathcal{W})$
for $i := \sqrt{-1}$, 
we can now consider the following 
modified weighted Bergman kernel 
$$
\beta_{q, \omega, \mathcal{W}} \;:= 
\; 2 \xi \, (1 + \frac{2q}{3}\,\Delta_{\omega})\, (B_{\xi ,\omega , \mathcal{Y}}
 - C_{\xi,\mathcal{Y}} ),
 \qquad \xi \in \Bbb C,
$$
where we used (5) and (6)
in defining $B_{\xi ,\omega , \mathcal{Y}}$ and $C_{\xi,\mathcal{Y}}$.
For the K\"ahler manifold $(M, \omega )$, consider the
Hamiltonian function $\hat{\sigma}_{\omega} \in C^{\infty}(M)_{\Bbb R}$ 
for $\mathcal{V}$ characterized by
$$
\mathcal{V}\,\lrcorner\,\omega \; =\; \bar{\partial} \hat{\sigma}_{\omega} 
\quad \text{ and }\quad
\int_M \sigma_{\omega} \, \omega^n \; =\; \int_M \hat{\sigma}_{\omega}  \,\omega^n.
$$
For each $\mathcal{Z}\in \frak z$, the associated Hamiltonian function $f^{}_{\mathcal{Z},\omega}\in C^{\infty}(M)_{\Bbb R}$ is uniquely characterized by 
 the identities $\mathcal{Z}\,\lrcorner\,\omega = \bar{\partial}  f^{}_{\mathcal{Z},\omega}$ and 
 $\int_M f^{}_{\mathcal{Z},\omega}\,\omega^n = 0$.
Note that $\mathcal{V} \in \frak z$. 
Then in this case also, as in (4), we obtain
\begin{equation}
\beta_{q,\omega, \mathcal{W}}\; 
=\; \sigma_{\omega} - \hat{\sigma}_{\omega} - f^{}_{\mathcal{W},\omega}
+ O(q). 
\end{equation}
Choose an arbitrary nonnegative integer $\ell$ with a real number $0<\alpha <1$.
For the space $\mathcal{F}_{\ell}$ of all $K$-invariant 
functions $f\in C^{\ell,\alpha}(M)_{\Bbb R}$ satisfying 
$\int_M f \omega_0^{\,n} = 0$, 
by setting $\mathcal{N}:= \operatorname{Ker} L_{\omega_0}\,\cap\,
\mathcal{F}_{\ell}$,
we have the identification
$$
\theta : \; \frak z \; \cong \; \mathcal{N},
\qquad \mathcal{Z} \,\leftrightarrow\, \theta (\mathcal{Z}) := f_{\mathcal{Z}, \omega_0},
$$
where $\mathcal{N}$ is independent of the choice of $\ell$.
Then the vector space $\mathcal{F}_{\ell}$ is written as a direct sum 
$\mathcal{N} \oplus \mathcal{N}_{\ell}^{\perp}$, where $\mathcal{N}_{\ell}^{\perp}$ 
is the space of all functions $f$ in $\mathcal{F}_{\ell}$ such that
 $\int_M  f \nu \omega_0^{\,n} = 0$ for all $\nu \in \mathcal{N}$.
 For an arbitrary integer $k  \geq 0$,
 we make a small perturbation of $\omega_0$ by varying $\omega$ in the space 
 $\{ \, \omega_0+\sqrt{-1}\,\partial\bar{\partial}\varphi \,;\, \varphi \in \mathcal{N}_{k+4}^{\perp}\,\}$. Since $\omega_0$ is an extremal K\"ahler metric, 
 we see from (8) that $\beta_{q,\omega, \mathcal{W}}$ 
 vanishes at  $(q, \omega , \mathcal{W}) = (0,\omega_0, 0)$, i.e., 
 $$
 \beta_{0, \omega_0, 0}\; =\; 0
 \qquad \text{ on $M$.}
 $$
Again from (8), we see that the Fr\'echet derivatives 
$D_\omega \beta_{q, \omega, \mathcal{W}}$ and
$D_{\mathcal{W}} \beta_{q, \omega, \mathcal{W}}$ 
of $\beta_{q, \omega, \mathcal{W}}$
at $(q, \omega , \mathcal{Y}) = (0,\omega_0, 0)$ are
$$\left\{\begin{array}{ll}
\;\;\;\left \{\,({D^{}_{\omega}\beta_{q,\omega, \mathcal{W}}^{}})\,
(\sqrt{-1}\,\partial\bar{\partial} \varphi )\,
\right\}_{| (q,\omega, \mathcal{W})= (0,\omega_0,0)}
\; &=\; L_{\omega_0} \varphi, \qquad  \varphi \in \mathcal{N}_{k+4}^{\perp}, \\
\quad\, ({D^{}_{\mathcal{W}}\beta_{q,\omega, \mathcal{W}}^{}})_{| (q,\omega, \mathcal{W}) 
= (0,\omega_0,0)}
\; &=\;\;\; -\,\theta, \qquad \quad \; \,\text{ on $\mathcal{N}$}.
\end{array}
\right. 
$$
Since $L_{\omega_0}: \mathcal{N}_{k+4}^{\perp} \to \mathcal{N}^{\perp}_k$ 
is invertible, and since $\theta$ is an isomorphism, the 
implicit function theorem is now applicable to the map: 
$(q, \omega, \mathcal{W}) \mapsto \beta_{q,\omega, \mathcal{W}}$. Then 
for some $0<\varepsilon \ll 1$, 
we can write
$$
\omega = \omega (q) \quad \text{ and }
\quad \mathcal{W} = \mathcal{W}(q),
\qquad 0\leq q \leq \varepsilon,
$$
solving the one-parameter family of equations
$$
\beta_{q,\omega, \mathcal{W}} \; =\; 0, 
\qquad q\geq 0,
$$
in $(\omega, \mathcal{W})$.  Hence by setting $\mathcal{Y}(q) := 
i q^2(\mathcal{V}+\mathcal{W}(q))$ for  $q = 1/m$, the $m$-th weighted Bergman kernel
$B_{m, \omega (q), \mathcal{Y}(q)}$ is 
constant on $M$ for all $m >1/\varepsilon$.
Then by Fact 4 in Section 5, $(M, L^m)$ is $T$-stable for all 
$m > 1/\varepsilon$.

\medskip\noindent
{\em Step 2}.  Though we assumed that $T$ (cf.~Section 5)
coincides with $Z^{\Bbb C}$ in the first step,
it is better to choose $T$ as small as possible.
Then for a sufficiently small positive real number $\varepsilon$,
 the algebraic torus $T$ in $Z^{\Bbb C}$ generated by
$$
\bigcup_{1/\varepsilon < m  \in \Bbb Z} 
\left \{\;\mathcal{V} + \mathcal{W}(q), \; i\, \mathcal{V} + \,i\,\mathcal{W}(q)\;\right \}
$$
is a good choice, where $q = 1/m$ and $i := \sqrt{-1}$.        

\section{Concluding remarks}

For Conjecture in Section 3,
``if'' part is known to be a difficult problem, and is of particular interest.
Assuming that a polarized manifold $(M, L)$ is asymptotically Chow-Mumford stable,
we are asked whether there exist K\"ahler metrics of constant scalar curvature in 
$c_1(L)_{\Bbb R}$. 
For simplicity, we consider the case where $H$ is trivial.
Then by Fact 3, the equation 
\begin{equation}
\beta_{q, \omega} = 0
\end{equation}
has solutions $(q, \omega ) = (1/m, \omega_m)$, $m \gg 1$, while  
each $\omega_m$ is an $m$-th balanced metric for $(M,L)$. 
Moreover, the triviality of $H$ 
implies that $\omega_m$ is the only $m$-th balanced metric for $(M,L)$. 
Then we are led to study the graph 
$$
\mathcal{E}\, :=\; \{\, (q, \omega )\,;\, \beta_{q,\omega} = 0\, \}
$$
in $\Bbb C \times \mathcal{K}$, where $\mathcal{K}$ is the space of
all K\"ahler metrics in the class $c_1(L)_{\Bbb R}$. 
Let $\mathcal{H}$, $\mathcal{M}$ be the sets of all Hermitian metrics for $L$, $V_m$, respectively.
Consider the 
Fr\'echet derivative $D_{\omega} \beta_{q,\omega}$ at $(q, \omega ) = (1/m, \omega_m)
\in \mathcal{E}$, where $m \gg 1$. 
With the same setting of differentiability as in Section 6, 
the Fr\'echet derivative will be shown to be invertible. 

\medskip
Let us give a rough idea how the invertibleness can be proved.
It suffices to show that
the operator ${D_{ \omega}B_{m,\omega}}$ is invertible at $\omega = \omega_m$.
Choose a Hermitian metric $h_m$ for $L$ such that 
$c_1(L; h_m ) = \omega_m$.
Put $ \Psi :=  \{\, \psi \in C^{\infty}(M)_{\Bbb R}\,;\, 
\int_M \psi \omega_m^{\;n} = 0\,\}$. Let $\varphi \in \Psi$. 
Then the $\Bbb R$-orbit in $\mathcal{H}$ through $h_m$ written in the form
\begin{equation}
  h_{\varphi ;t} := \,e^{-t\varphi}h_m, \qquad t\in \Bbb R,
\end{equation}
projects to the $\Bbb R$-orbit
$ \omega_{\varphi,t} := \omega_m + \sqrt{-1}\, t \,\partial\bar{\partial}\varphi$, $t \in \Bbb R$,
in $\mathcal{K}$ through $\omega_m$.
Note also that every Hermitian metric $h$ for $L$ induces a Hermitian pairing $<\; ,\;>^{}_{L^2}$
which will be denoted by $(V_m, \tilde{h})$ (cf. Section 4).
Choose an orthonormal basis $\{\sigma_1,\sigma_2,\dots,\sigma_{N_m}\}$ for $(V_m, \tilde{h}_m)$. Let $\Phi$ be the space of all $\varphi \in \Psi$ such that  
$$
(d/dt)(\tilde{h}_{\varphi,t})_{|t=0} \; =\; 0
\qquad \text{ on $V_m$.}
$$ 
In other words, if $\varphi \in \Phi$, then
the basis $\{\sigma_1,\sigma_2,\dots,\sigma_{N_m}\}$ 
infinitesimally remains to be an orthonormal basis for $(V_m, \tilde{h}_{\varphi,t})$ 
at $t  =0$ with $t$ perturbed a little.
Since
$$
B_{m,\omega} := |\sigma_1|_h^2+|\sigma_2|_h^2 + \dots + |\sigma_{N_m}|_h^2
$$
is obtained from the contraction of $\Sigma:= \sigma_1\bar{\sigma}_1 + \sigma_2\bar{\sigma}_2 + \dots + 
\sigma_{N_m}\bar{\sigma}_{N_m}$ by $h^m$, 
and since $\Sigma$ is fixed  
by the infinitesimal action $(d/dt)(h_{\varphi,t})$ at $t = 0$ in (10), 
we obtain
\begin{equation}
(d/dt)(B_{m, \omega^{}_{\varphi,t}})_{|t=0}\; =\; (d/dt)(h^m_{\varphi,t})_{|t=0} \;=\; - m\varphi.
\end{equation}
Let $\Bbb H$ denote the set of all Hermitian metrics on $V_m$. Then we have a natural
projection $\pi : \Psi \to T_{\tilde{h}_m}\Bbb H$  defined by 
$$
\qquad (d/dt)(h_{\varphi,t})_{|t=0}\;  \mapsto \; (d/dt)(\tilde{h}_{\varphi,t})_{|t=0}.
$$
In view of $\dim T_{\tilde{h}_m}\Bbb H < + \infty$,  it is now easy to see that this map is surjective.
Since the kernel of this map is $\Phi$, we obtain 
\begin{equation}
\Psi/\Phi \; \cong \; T_{\tilde{h}_m}\Bbb H.
\end{equation}
Now by (11) and (12), 
the uniquness of a balanced metric (by $H= \{1\}$)
implies the required invertibleness of the Fr\'echet derivative.
(see also Donaldson \cite{Don4}).
\qed

\medskip
Now we can apply the implicit function theorem
to obtain an open neighborhood 
$U$ of $1/m$ in $\Bbb C$ such that, 
for each $\xi \in U$, there exists a unique K\"ahler metric $\omega (\xi )$
in $c_1(L)_{\Bbb R}$ satisfying
$\beta_{\xi,\omega (\xi )} = 0$.

\medskip
Assuming nonexistence of K\"ahler metrics of constant scalar curvature, 
we have some possibility that, by
an argument as in Nadel, destabilizing objects are obtained
by studying the behavior of the solutions along the boundary point of $\mathcal{E}$.

\medskip
Finally, there are many interesting topics, which are related to ours,
such as the geometry of K\"ahler potentials \cite{Che}, \cite{CC} \cite{PS1} \cite{Sem} (see also \cite{Mab1}).
The uniqueness of extremal K\"ahler metrics modulo the action of holomorphic automorphisms in compact cases
are recently given by \cite{CT2} (cf.~\cite{BM},~\cite{Don2}).
New obstructions to semistability of manifolds or to the existence of extremal metrics
are done by \cite{Fut2}, \cite{Mab2}, \cite{CT2}. 
The concept of K-stability, introduced by Tian \cite{Tia3} and reformulated by Donaldson \cite{Don3},
is deeply related to the Hilbert-Mumford stability criterion, and various kinds of works 
are actively done related to algebraic geometry.

\end{document}